\documentclass[12pt]{article}

\usepackage{amsmath}
\usepackage{amsfonts}
\usepackage{amssymb}
\usepackage{graphicx, color}
\usepackage{hyperref}
\usepackage{epsfig}
\usepackage{epstopdf} 
\usepackage{esint}

\definecolor{darkgreen}{rgb}{0,0.55,0}

\newtheorem{theorem}{Theorem}[section]
\newtheorem{lemma}[theorem]{Lemma}
\newtheorem{corollary}[theorem]{Corollary}
\newtheorem{remark}[theorem]{Remark}

\newtheorem{definition}{Definition}

\DeclareSymbolFont{AMSb}{U}{msb}{m}{n}
\DeclareMathSymbol{\N}{\mathbin}{AMSb}{"4E}
\DeclareMathSymbol{\Z}{\mathbin}{AMSb}{"5A}
\DeclareMathSymbol{\R}{\mathbin}{AMSb}{"52}
\DeclareMathSymbol{\Q}{\mathbin}{AMSb}{"51}
\DeclareMathSymbol{\I}{\mathbin}{AMSb}{"49}

\newcommand{\calA}{{\mathcal A}}
\newcommand{\calH}{{\mathcal H}}

\begin{document}

\title{Symmetry of solutions of a mean field equation on flat tori }
\author{{Changfeng Gui\footnote{Department of Mathematics, University of Texas at San Antonio, Texas, USA. E-mail: changfeng.gui@utsa.edu.  }
\qquad Amir Moradifam\footnote{Department of Mathematics, University of California, Riverside, California, USA. E-mail: moradifam@math.ucr.edu. }}}
\date{\today}

\smallbreak \maketitle

\begin{abstract}
We study symmetry of solutions of the mean field equation 
\[ \Delta u +\rho(\frac{Ke^u}{\int_{T_\epsilon}  Ke^u} -\frac{1}{|T_\epsilon|} )=0\]
on the flat torus $T_\epsilon=[-\frac{1}{2\epsilon}, \frac{1}{2\epsilon}] \times [-\frac{1}{2}, \frac{1}{2}]$ with $0<\epsilon \leq 1$, where $K\in C^2({T}_\epsilon)$ is a positive function with $-\Delta \ln K \leq \frac{\rho}{|T_\epsilon|}$  and $\rho \leq 8\pi$.  We prove that if $(x_0,y_0)$ is a critical point of the function $u+ln(K)$, then $u$ is evenly symmetric about the lines $x=x_0$ and $y=y_0$, provided $K$ is evenly symmetric about these lines. In particular we show that all solutions  are one-dimensional if $K\equiv 1$ and $\rho \leq 8\pi$. The results are sharp and answer a conjecture of Lin and Lucia affirmatively. We also prove some symmetry results for mean field equations on annulus. 

\end{abstract}
\maketitle

\section{Introduction} 

Consider the mean field equation 
 
 \begin{equation}\label{meanfield0}
 \Delta u +\rho(\frac{e^u}{\int_{T_\epsilon}  e^u} -\frac{1}{|T_\epsilon|} )=0, \quad (x, y) \in T_\epsilon,
\end{equation}
on the flat torus with fundamental domain 
\begin{equation}\label{torus}
T_\epsilon=[-\frac{1}{2\epsilon}, \frac{1}{2\epsilon}] \times [-\frac{1}{2}, \frac{1}{2}] , \ \ 0<\epsilon \leq 1.
\end{equation}
Define 
\[\calH(T_\epsilon)=\{u\in H^1(T_\epsilon): \ \ \int_{T_\epsilon}u =0\}.\]
Solutions of (\ref{meanfield0}) are critical points of the functional $F_{\rho, \epsilon}: \calH(T_\epsilon)\rightarrow \R$ defined by 
\[F_{\rho, \epsilon}:=\frac{1}{2}\int_{T_\epsilon}|\nabla u|^2-\rho \ln(\frac{1}{|T_\epsilon|}\int_{T_\epsilon}e^u).\]
Notice that, since both equation (\ref{meanfield0}) and  functional $F_{\rho,\epsilon}$ are invariant under dilations, we do not lose generality by taking the vertical size of $T_\epsilon$ to be $1$. 

Equations of type  (\ref{meanfield0}) arise in Onsager's vortex theory for one specie 
(see,  \cite{CLMP1-MR1145596,CLMP2-MR1362165,CK-MR1262195,K-MR1193342}). It is also obtained in the context of Chern-Simons guage theory (see \cite{Caff-Yang-MR1324400, Dunne, Hong-MR1050529,Jackiw-MR1050530, Tarantello-MR1400816}, etc). Tarantello \cite{Tarantello-MR1400816} showed that when the Chern-Simons coupling constant tends to zero, the asymptotic behavior of a class of solutions is described by 
\begin{equation}\label{meanfield}
 \Delta u +\rho(\frac{Ke^u}{\int_{T_\epsilon}  Ke^u} -\frac{1}{|T_\epsilon|} )=0, \quad (x, y) \in T_\epsilon,
\end{equation}
where $\rho=4\pi N$, $N$ is an integer called the vortex number, and $K$ is a prescribed non-negative function. 

Various results have been obtained regarding the existence and qualitative properties of mean field equations on flat tori (see, e.g., \cite{Cabre-MR2180306, Chen-Lin-MR1831657, Ding1-MR1491984,Ding3-MR1624438,Ding2-MR1712560, K-MR1193342, YLi-MR1673972, LS-MR1322618,  LL1-MR2265623, LL2-MR2352519, N-Tarantello-MR1664542, Ricciardi-Tarantello-MR1664758, Struwe-Tarantello-MR1619043}). If $\rho \leq 0$, it is easy to see that the functional $F_{\rho,\epsilon}$ is strictly convex and consequently zero is the only solution of (\ref{meanfield0}) which is also a minimizer of $F_{\rho,\epsilon}$. For $ \rho>0 $ it follows from the Moser-Trudinger inequality \cite{Moser-MR0301504, T-MR0216286} that the functional $F_{\rho,\epsilon}$ is bounded from below if and only if $\rho \leq 8\pi$. Moreover if $\rho <8 \pi$, then $F_{\rho,\epsilon}$ is coercive and therefore it admits a minimizer. For $\rho= 8\pi$ existence of a minimizer is discussed in \cite{Ding1-MR1491984} and \cite{N-Tarantello-MR1664542}. 

For $\rho > 8\pi$, the functional $F_{\rho, \epsilon}$ is unbounded below. It is shown in \cite{Ding2-MR1712560} by a minmax method that  \eqref{meanfield}  has a solution for $ \rho \in (8\pi, 16 \pi)$.  In \cite{Struwe-Tarantello-MR1619043} it is  proved that if $\epsilon=1$  and $8\pi <\rho < 4\pi^2$, then the equation (\ref{meanfield0}) has two-dimensional solutions and $u\equiv 0$ is only a local minimizer of $F_{\rho,1}$. On the other hand it is proved in \cite{Ricciardi-Tarantello-MR1664758} that $\rho >4\pi^2 \epsilon$ is a necessary and sufficient condition for the existence of at least one non-zero one-dimensional solution. Moreover if $4\pi^2 \epsilon< 8\pi$ and $\rho \in (4\pi^2 \epsilon, 8\pi]$, then any minimizer is nonzero \cite{Ricciardi-Tarantello-MR1664758}. 

As stated above, there exist two-dimensional solutions for $\rho \in (8\pi, 4\pi^2)$ and $\epsilon=1$. So it is a natural question to ask whether all solutions of (\ref{meanfield0}) are one-dimensional for $\rho \leq 8\pi$. 
In \cite{Cabre-MR2180306}, Cabr\'{e}, Lucia, and Sanch\'{o}n proved that if \\ \\
\[\rho \leq \rho^*:=\frac{16 \pi^3}{\pi^2+\frac{2}{R_\epsilon^2}+\sqrt{(\pi^2+\frac{2}{R^2_{\epsilon}})^2-\frac{8\pi^3}{|T_\epsilon|}}} \leq 0.879 \times 8\pi,\] \\ \\
then every solution $u$ of (\ref{meanfield0}) depends only on $x$-variable. Here $R_\epsilon$ is the maximum conformal radius of the rectangle $T_\epsilon$.  

In \cite{LL1-MR2265623},  Lin and Lucia  proved that  constants are the only solutions of (\ref{meanfield0}) whenever \\ 
\begin{eqnarray*}
\rho \leq \left\{ \begin{array}{ll}
 8\pi &\text{if } \epsilon \geq \frac{\pi}{4}\\
32\epsilon &\text{if } \epsilon \leq \frac{\pi}{4}
\end{array} \right.
\end{eqnarray*}
\\ (see Theorem 1.3 in \cite{LL1-MR2265623}). This result is optimal only if $\epsilon \geq \frac{\pi}{4}$. 

Later in \cite{LL2-MR2352519} Lin and Lucia obtained optimal symmetry results for minimizers of the functional $F_{\rho, \epsilon}$. Indeed they proved the following theorem. \\ \\
{\bf Theorem A.} (Theorem 1.2 in \cite{LL2-MR2352519}) Let $T_\epsilon$ be the flat torus defined in (\ref{torus}) and suppose $\rho \leq 8\pi $. Then any global minimizer of $F_{\rho,\epsilon}$ is one-dimensional. In addition for $\rho \leq \min \{8\pi, 4\pi ^2 \epsilon\}$, $u\equiv 0$ is the unique golobal minimizer of the functional $F_{\rho,\epsilon}$. \\ \\
However, one-dimensional symmetry of solutions of (\ref{meanfield0}) still remained  open. In particular, Lin and Lucia \cite{LL2-MR2352519} conjectured that $u\equiv 0$ is the unique solution of (\ref{meanfield0}) whenever $\rho \leq \min \{8\pi, 4\pi ^2 \epsilon\}$. In this paper, among other results, we prove this conjecture. Indeed we prove the following optimal result which improves the results in \cite{Cabre-MR2180306}, \cite{LL1-MR2265623}, and \cite{LL2-MR2352519}, discussed above. \\

\begin{theorem}\label{meanFeildThoremTorus}
Suppose $\rho \leq 8\pi$ and let $u$ be a  solution of  \eqref{meanfield0}.  Then $u$ must be  one-dimensional.  In particular,  $u$ is constant if $ \rho \le \min\{ 8 \pi,  4 \pi^2 \epsilon\}. $
 \end{theorem}
 
We also prove the following theorem for the general mean field equation (\ref{meanfield}). \\
\begin{theorem}\label{meanFeildTehoremTorusGeneral}   Let $\rho \leq 8\pi$, $T_\epsilon$ be a flat torus defined in (\ref{torus}), and $K\in C^2(T_\epsilon)$ be a positive function with $-\Delta \ln K \leq \frac{\rho}{|T_\epsilon|}$ in $T_\epsilon$. Suppose that  $K$ is evenly symmetric in $x, y$, and $u$ is a solution of  \eqref{meanfield} with the  origin being  a critical point of $u$. Then $u$ is evenly symmetric  in $x$ and $y$, i.e., 
\[ u(x, y)=u(-x, y)=u(x, -y) ,  \ \ \forall (x, y) \in T_\epsilon.  \]
\end{theorem}
\vspace{0.4cm}

Letting $v=\ln K +u$, we obtain the following corollary of Theorem \ref{meanFeildTehoremTorusGeneral}. \\

\begin{corollary}\label{meanFeildTehoremTorusGeneral1}   Let $\rho \leq 8\pi$, $T_\epsilon$ be a flat torus defined in (\ref{torus}), and $K\in C^2(T_\epsilon)$ be a positive function with $-\Delta \ln K =C\leq \frac{\rho}{|T_\epsilon|}$ in $T_\epsilon$, where $C$ is a constant.  Suppose $(x_0,y_0)$ is a critical point of $v=\ln K+u$. Then  $\ln K+u$ is symmetric with respect to the lines $x=x_0$ and $y=y_0$. 
\end{corollary}
\vspace{0.3cm}

We also present symmetry results for mean field equations on annulus in Section 3, improving corresponding results in \cite{Chen-Lin-MR1831657}. 

\section{Proof of the symmetry results on tori}
In this section we present the proofs of our main results, Theorems \ref{meanFeildThoremTorus} and \ref{meanFeildTehoremTorusGeneral}. The proofs are based on the Sphere Covering Inequality recently proved by the authors in \cite{GM-SphereCovering}. \\ \\
{\bf Theorem B.} (Theorem 3.1 in \cite{GM-SphereCovering})
Let $\Omega $ be a simply-connected subset of $R^2$ and assume $w_i \in C^2(\overline{\Omega})$, $i=1,2$ satisfy
\begin{equation}
\Delta w_i +e^{w_i}=f_{i},
\end{equation}
where $f_2 \geq 0$  and $f_2 \geq f_1 $ in $\Omega$. If $w_2\not\equiv w_1$ in $\Omega$ and $w_2=w_1$ on $\partial \Omega$, then
\begin{equation}
\int_{\Omega} (e^{w_1}+e^{w_2}) \geq 8\pi. 
\end{equation}
Moreover,  the equality only holds when $f_2 \equiv f_1 \equiv 0$ and $(\Omega, e^{w_i}dy ) $, $i=1,2$ are isometric to  two complimenting spherical caps on the standard sphere with radius $\sqrt2$.  \\ \\ 
{\bf Proof of Theorem \ref{meanFeildTehoremTorusGeneral}.} Let $v:=u+\ln (K)+\ln \rho -\ln (\int_{T_\epsilon} K e^u)$. Then $v$ satisfies 
\begin{equation}\label{gaussian}
\Delta v+e^v=\frac{\rho}{|T_\epsilon|}+\Delta \ln K \geq 0 \quad (x, y) \in T_\epsilon.
\end{equation}
We claim that  $u$ is even in $x$ and $y$, i.e., 
\[ v(x, y)=v(-x, y)=v(x, -y) ,  \ \ \forall (x, y) \in T_\epsilon.\]
To prove the claim, define 
\[w(x, y)= v(x, y)-v(x, -y),  \ \  (x, y) \in T_\epsilon.\]
Suppose that $w \not \equiv 0$,  and set 
\[\Omega^+:=\{x\in  [-\frac{1}{2\epsilon},  \frac{1}{2\epsilon}] \times [0, \frac{1}{2}] : \ \ w(x)>0 \}\]
and
\[\Omega^-:=\{x\in  [-\frac{1}{2\epsilon},  \frac{1}{2\epsilon}] \times [0, \frac{1}{2}] : \ \ w(x)<0 \}.\]
Note that $w=0$ on $\Gamma_0 \cup \Gamma_1$, where 
\[\Gamma_0=\{(x, 0): -\frac{1}{2\epsilon}  \le x \le \frac{1}{2\epsilon}\} \ \ \hbox{and} \ \   \Gamma_1=\{(x, 1): -\frac{1}{2\epsilon}  \le x \le \frac{1}{2\epsilon}\}.\]
Since $(0,0)$ is a critical point of $u$ and $K$ is evenly symmetric in $y$,  
\[\frac{\partial}{\partial y} w(0,0)=0.\]
Therefore, by Hopf's lemma, the nodal set of $w$ must contain a curve $\Gamma$ originating from the origin and having a transversal intersection with  $\Gamma_0$.   We now discuss two cases:\\ \\
(a)  $\Gamma $ reaches the  boundary curve $\Gamma_1$;   \\ 
(b)  $\Gamma$ does not reach the boundary curve $\Gamma_1$.\\ 

In case (a),  there must be another nodal curve of $w$ connecting  $\Gamma_1$ and $\Gamma_2$,  and hence  there are at least two simply connected 
nonempty regions $\Omega_1, \Omega_2 \subset [-\frac{1}{2\epsilon},  \frac{1}{2\epsilon}] \times [0, \frac{1}{2}]  $ such that $\Omega_1 \cap \Omega_2=\emptyset$ and $w=0 $ on $\partial \Omega_1 \cup \partial \Omega_2$. Hence on each $\Omega_i$, $i=1,2$,  equation (\ref{gaussian}) has two distinct solutions $v(x,y)$ and $v(x,-y)$ with $v(x,y)=v(x,-y)=0$ on $\partial \Omega_i$. 
Therefore,  by the Sphere Covering Inequality (Theorem B)  we conclude that 
$$
\int_{\Omega_i} (e^{v(x,y)}+e^{v(x, -y)} )> 8 \pi, \quad i=1, 2.
$$
Hence 
\begin{equation}\label{16pi}
\rho =\int_{T_\epsilon}  e^v \ge  \sum_{i=1}^{2} \int_{\Omega_i} (e^{v(x,y)}+e^{v(x, -y)} )> 16 \pi,
\end{equation}
which is contradiction. \\ 

In Case (b),  there exists at least one simply-connected region $\Omega_1 \subset \subset [-\frac{1}{2\epsilon},  \frac{1}{2\epsilon}] \times [0, \frac{1}{2}] $  such that $w=0 $ on $\partial \Omega_1$.  Thus the equation (\ref{gaussian}) has two distinct solutions $v(x,y)$ and $v(x,-y)$ on $\Omega_1$ with $v(x,y)=v(x,-y)=0$ on $\partial \Omega_1$. 
Therefore,  by the Sphere Covering Inequality (Theorem B) again  we conclude that 
$$
\int_{\Omega_1} (e^{v(x,y)}+e^{v(x, -y)} )> 8 \pi.
$$
Consequently 
\begin{equation}\label{8pi}
\rho =\int_{T_\epsilon}  e^v \ge  \int_{\omega_1} (e^{v(x,y)}+e^{v(x, -y)} )> 8 \pi,
\end{equation}
which is a contradiction. 
In both cases,  we conclude $\rho >8\pi$ which contradicts the assumption $\rho \leq 8\pi$.  Hence  $u(x, -y)\equiv u(x, y)$ in $T_\epsilon$.  The proof of $u(x, y)\equiv u(-x, y)$ in $T_\epsilon$ is similar. \hfill $\Box$ 

\begin{remark}
Notice that if $u$ has a critical point at the origin and another critical point $X^*=(x^*,1)$ with $-\frac{1}{2\epsilon}\leq x^* \leq \frac{1}{2\epsilon}$, then case (b) in the proof of Theorem \ref{meanFeildTehoremTorusGeneral} can not happen and therefore $u$ must be symmetric with respect to $x$-axis if $\rho \leq 16\pi $. Similarly if $u$ has another critical point $Y^*=(\frac{1}{2\epsilon},y^*)$ with $-\frac{1}{2}\leq y^* \leq \frac{1}{2}$, then $u$ will be symmetric with respect to $y$-axis for $\rho \leq 16 \pi$. In particular, if $u$ has a critical point at the origin and another critical point at the corner of the torus, then $u$ must be symmetric with about $x$ and $y$-axis for $\rho \leq 16 \pi$. \\ 
\end{remark}
From the above proof, we can easily see the following. \\
\begin{remark}
Let $\rho \leq 8\pi$, $T_\epsilon$ be a flat torus defined in (\ref{torus}), and $K\in C^2(T_\epsilon)$ be a positive function with $-\Delta \ln K \leq \frac{\rho}{|T_\epsilon|}$ in $T_\epsilon$. If $K$ is only evenly symmetric in $y$ and $u_y(x^*, 0)=0$ for some $x^* \in [-\frac{1}{2\epsilon},  \frac{1}{2\epsilon}]$,  then $u$ is evenly symmetric in $y$.  Similarly, if $K$ is only evenly symmetric in $x$ and $u_x(0, y^*)=0$ for some $y^* \in [-\frac{1}{2},  \frac{1}{2}]$, then $u$ is evenly symmetric in $x$. \\ 
\end{remark}
Next we prove Theorem \ref{meanFeildThoremTorus}.  Following \cite{LL2-MR2352519} let us first define Steiner symmetric solutions. \\ 
\begin{definition}
Let $T=(-a,a)\times (-b,b)$ be a flat torus. A function $u \in H^1(T)$ is said to be Steiner symmetric on $T$ if 

\begin{eqnarray*}
\left\{ \begin{array}{lll}
u(x,y)=u(-x,y)=u(x,-y) & \forall (x,y)\in T, \\
\frac{\partial u}{\partial x}(x)\leq 0 & \forall (x,y) \in (0,a) \times (-b,b),\\
\frac{\partial u}{\partial y}(x)\leq 0 & \forall (x,y) \in (-a,a) \times (0,b). \\ \\
\end{array} \right. \\ 
\end{eqnarray*} 
\end{definition}
Now we are ready to prove Theorem \ref{meanFeildThoremTorus}. \\ \\
{\bf Proof of Theorem \ref{meanFeildThoremTorus}.} Without loss of generality we may assume that $u$ has a maximum point at $(0,0)$. Thus it follows from Theorem \ref{meanFeildTehoremTorusGeneral} that $u$ is symmetric about $x$ and $y$-axis. 

Assume now that $u_{x}(x^*, y^*)=0 $ for some $x^* \in (0,  \frac{1}{2\epsilon}) $. Reflecting $u$ about $x=x^*$ and applying an argument similar to the one in the proof of  Theorem \ref{meanFeildTehoremTorusGeneral}, we can conclude that $u$ is evenly symmetric about $x=x^*$.  Similarly if $u_{y}(x^*, y^*)=0 $, then we conclude that $u$ is symmetric about $y=y^*$. Therefore, there exist positive integers $l, m$ such that  $u$ is a periodic function with periods $\frac{1}{l\epsilon}, \frac{1}{m}$ in $x, y$ variables, respectively. Moreover 
\[u_{x}(x, y)\leq 0,  u_{y}(x, y) \leq 0  \ \ \hbox{ for}\ \   (x, y) \in  [0, \frac{1}{2l\epsilon}] \times [0, \frac{1}{2m}].\]
In other words,  $u$ can be regarded as a Steiner symmetric solution on a flat 
torus with fundamental domain as $[-\frac{1}{2l\epsilon}, \frac{1}{2l\epsilon}] \times [-\frac{1}{2m}, \frac{1}{2m}]$. Since $\rho \leq 8\pi$, it follows from  Theorem 1.2  b) in  \cite{LL2-MR2352519}  that $u$ is one-dimensional
in $[-\frac{1}{2l\epsilon}, \frac{1}{2l\epsilon}] \times [-\frac{1}{2m}, \frac{1}{2m}]$, and consequently $u$ must be one-dimensional in  $T_\epsilon$. 

Riccardi and Tarantello in \cite{Ricciardi-Tarantello-MR1664758} showed that  $\rho >4 \pi^2 \epsilon$ is a necessary and sufficient condition for the existence of a non constant one dimensional solution. Hence  $u$ must be constant if $ \rho \le \min\{ 8 \pi,  4 \pi^2 \epsilon\}$.  The proof is now complete. \hfill $\Box$ \\  

\section{Mean field equations on annulus}
In this section we prove symmetry results for mean field equations on an annulus. Let $\mathcal{A}$ be an arbitrary annulus in $\R^2$, i.e.  
\begin{equation*}
\mathcal{A}:= \{(x,y) \in \R^2: \ \  a < |(x, y)|< b \} \ \ \ \ \hbox{for some}\ \ a<b. 
\end{equation*}
We consider the following mean field equation 

\begin{eqnarray}\label{AnnulusPDE}
\left\{ \begin{array}{ll}
 \Delta u+\rho \frac{K(x,y)e^u}{\int_{\mathcal{A}}K(x,y)e^u dx}=f(x,y) \geq 0 &\text{in} \ \ \mathcal{A}\\
u(x,y)=\alpha &\text{if } |(x, y)|=a\\
u(x,y)=\beta &\text{if } |(x, y)|= b,
\end{array} \right.
\end{eqnarray}
where $K $ is a  positive radial function with $\Delta \ln K \ge 0$,  $f$ is a non-negative radial function,  and  $\alpha ,\beta \in \R$.

For $\rho \in (8\pi, 16\pi)$,  it is shown in \cite{Ding2-MR1712560} that there exists a solution to \eqref{AnnulusPDE} when $\alpha=\beta=0$.  In  \cite{Chen-Lin-MR1831657}, it is proved that if $\rho_i \to 8\pi$ from left and $\rho_i \to 16 \pi$ from right, then for large $i$ the solutions have  even symmetry about a line passing through the origin,  when blow-up of the solution happens.  Below  (see Theorem \ref{AnnulusTheorem}) we shall show  that the solutions always have  even symmetry. \\

We will need the following lemma (see, e.g. \cite{BL-MR3201892, GL-MR2670931, LL2-MR2352519}). 

\begin{lemma}\label{4piBound}
Let $\Omega \subset \R^2 $ be a simply-connected domain and  assume that $w\in C^2(\overline{\Omega})$ satisfies 
$\Delta w+e^{w} \ge 0$ in $\overline{\Omega}$ and $\int_{\Omega} e^w \le 8 \pi$.   Consider an open set $\omega \subset \Omega $  and define
the first eigenvalue of the operator $ \Delta + e^w $ in $H_0^1(\omega) $ by
$$ 
\lambda_{1, w} (\omega):=  \inf_{\phi \in H^1_0 (\omega) }  \bigl(\int_{\omega} |\nabla \phi| ^2 - \int_{\omega} \phi^2 e^w \bigr)   \le 0.
$$
Then $\int_{\Omega}e^{w} \ge 4 \pi$  if $\lambda_{1, w} (\Omega) \le 0$.\\
\end{lemma}

\begin{theorem} \label{AnnulusTheorem}
Assume $\rho < 16 \pi$. Then, afert a proper rotation, every solution of (\ref{AnnulusPDE}) must be evenly symmetric about $x$-axis.  Furthermore,  either  $u$ is radially symmetric or the angular derivative $u_{\theta}$ of $u$  doesn't change sign in
\begin{equation*}
\mathcal{A^+}:= \{(x, y) \in \R^2: \ \  a < |(x, y)|< b, \quad  y>0 \}. \\
\end{equation*}
\end{theorem}
{\bf Proof.} Define 
\[v=u+ln(K)+\ln(\rho)-\ln(\int_{\mathcal{A}}Ke^u dx).\]
Then $v$ satisfies 
\begin{equation}\label{Annulus-v-PDE}
\Delta v+e^{v}=f+\Delta \ln K \geq 0.
\end{equation}
Fix $x_0 \in (a, b)$.  Without loss of generality,   we may assume that at $u$ attains its minimun  on the circle $|(x, y)|=x_0$ at the point $(x_0, 0)$. Note that 
 $\frac{\partial v}{\partial \theta}(x_0, 0)=0$. Define 
\[w(x,y):=v(x,y)-v(x,-y).\] 
Since 
\[\frac{\partial v}{\partial \theta}(x_0,0)=\frac{\partial v}{\partial y}(x_0,0)=0,\]
it follows from the Hopf's lemma that the nodal line of $w$ divided a neighborhood of $x_0$ into at least four regions. Hence there exists two simply connected regions 
\[\Omega_1, \Omega_2 \subset \{(x,y)\in \mathcal{A}: y>0\}\]
such that $w=0$ on $\partial \Omega_1 \cup \partial \Omega_2$. Therefore on each $\Omega_i$, $i=1,2$, the equation (\ref{Annulus-v-PDE}) has two solutions $v(x,y)$ and $v(x,-y)$ with  $v(x,y)=v(x,-y)$ on $\partial \Omega_i$. Thus it follows from the Sphere Covering Inequality (Theorem B) that 
\[\rho =\int_{\mathcal{A}}  e^v \ge  \sum_{i=1}^{2} \int_{\Omega_i} (e^{v(x,y)}+e^{v(x, -y)} )\geq 16 \pi.\]
This is a contradiction to the assumption $\rho< 16 \pi$.  Therefore, we conclude that $w \equiv 0$ in $\calA$, 
and $u$ is evenly symmetric about $x$-axis.  Furthermore,  $\phi:=v_\theta=u_\theta$ satisfies the linearized equaiton
\begin{eqnarray}\label{AnnulusLinearPDE}
\left\{ \begin{array}{ll}
 \Delta \phi+ e^v \phi = 0 &\text{in} \ \ \mathcal{A}\\
\phi(x, y) =0 &\text{if } |(x, y)|=a,  \text{ or } |(x, y)|= b, \text{ or }  y=0. 
\end{array} \right.
\end{eqnarray}
Assume $u$ is not radially symmetric.  If $\phi$ changes sign
in $\calA^{+}$,  then there are at least two regions $\Omega_1, \Omega_2 \subset \calA^+$ with $\phi=0$ on $\partial \Omega_1 \cup \partial \Omega_2$. Hence it follows from Lemma \ref{4piBound} that 
\[\rho =2\int_{\mathcal{A}^+}  e^v \ge  2 \sum_{i=1}^{2} \int_{\Omega_i} e^{v}\geq 16 \pi.\]
This is a contradiction. Hence $u_\theta$ does not change sign in $\calA^+$. 
The proof is now complete. \hfill $\Box$
\\ 

\begin{remark}
If $f+\Delta \ln K \not \equiv 0$ in $\mathcal{A}$, then the condition $\rho <16 \pi$ in the statement of Theorem \ref{AnnulusTheorem} can be replaced by $\rho \leq 16 \pi$. 
\end{remark}

\begin{remark}
Assume that  $K, f$ are only evenly symmetric about $x$-axis in Theorem \ref{AnnulusTheorem}. 
The above proof also indicates that $u$ must be evenly symmetric about $x$-axis if there is a point $(x_0, 0)$ on $x$-axis such that $u_y(x_0, 0)=0$. 
\end{remark}

\bibliographystyle{plain}
\bibliography{FlatTorus}

\begin{thebibliography}{10}

\bibitem{BL-MR3201892}
Daniele Bartolucci and Chang-Shou Lin.
\newblock Existence and uniqueness for mean field equations on multiply
  connected domains at the critical parameter.
\newblock {\em Math. Ann.}, 359(1-2):1--44, 2014.

\bibitem{Cabre-MR2180306}
Xavier Cabr{\'e}, Marcello Lucia, and Manel Sanch{\'o}n.
\newblock A mean field equation on a torus: one-dimensional symmetry of
  solutions.
\newblock {\em Comm. Partial Differential Equations}, 30(7-9):1315--1330, 2005.

\bibitem{Caff-Yang-MR1324400}
Luis~A. Caffarelli and Yi~Song Yang.
\newblock Vortex condensation in the {C}hern-{S}imons {H}iggs model: an
  existence theorem.
\newblock {\em Comm. Math. Phys.}, 168(2):321--336, 1995.

\bibitem{CLMP1-MR1145596}
E.~Caglioti, P.-L. Lions, C.~Marchioro, and M.~Pulvirenti.
\newblock A special class of stationary flows for two-dimensional {E}uler
  equations: a statistical mechanics description.
\newblock {\em Comm. Math. Phys.}, 143(3):501--525, 1992.

\bibitem{CLMP2-MR1362165}
E.~Caglioti, P.-L. Lions, C.~Marchioro, and M.~Pulvirenti.
\newblock A special class of stationary flows for two-dimensional {E}uler
  equations: a statistical mechanics description. {II}.
\newblock {\em Comm. Math. Phys.}, 174(2):229--260, 1995.

\bibitem{CK-MR1262195}
Sagun Chanillo and Michael Kiessling.
\newblock Rotational symmetry of solutions of some nonlinear problems in
  statistical mechanics and in geometry.
\newblock {\em Comm. Math. Phys.}, 160(2):217--238, 1994.

\bibitem{Chen-Lin-MR1831657}
Chuin~Chuan Chen and Chang-Shou Lin.
\newblock On the symmetry of blowup solutions to a mean field equation.
\newblock {\em Ann. Inst. H. Poincar\'e Anal. Non Lin\'eaire}, 18(3):271--296,
  2001.

\bibitem{Ding1-MR1491984}
Weiyue Ding, J{\"u}rgen Jost, Jiayu Li, and Guofang Wang.
\newblock The differential equation {$\Delta u=8\pi-8\pi he^u$} on a compact
  {R}iemann surface.
\newblock {\em Asian J. Math.}, 1(2):230--248, 1997.

\bibitem{Ding3-MR1624438}
Weiyue Ding, J{\"u}rgen Jost, Jiayu Li, and Guofang Wang.
\newblock An analysis of the two-vortex case in the {C}hern-{S}imons {H}iggs
  model.
\newblock {\em Calc. Var. Partial Differential Equations}, 7(1):87--97, 1998.

\bibitem{Ding2-MR1712560}
Weiyue Ding, J{\"u}rgen Jost, Jiayu Li, and Guofang Wang.
\newblock Existence results for mean field equations.
\newblock {\em Ann. Inst. H. Poincar\'e Anal. Non Lin\'eaire}, 16(5):653--666,
  1999.

\bibitem{Dunne}
G.~Dunne.
\newblock Self-dual chern simons theories.
\newblock {\em Lecture Notes in Physics. New Series M 36. New York:
  Springer-Verlag}, 1996.

\bibitem{GL-MR2670931}
Nassif Ghoussoub and Chang-Shou Lin.
\newblock On the best constant in the {M}oser-{O}nofri-{A}ubin inequality.
\newblock {\em Comm. Math. Phys.}, 298(3):869--878, 2010.

\bibitem{GM-SphereCovering}
Changfeng Gui and Amir Moradifam.
\newblock The sphere covering inequality and its applications.
\newblock {\em http://arxiv.org/pdf/1605.06481v1.pdf, Submitted}, 2016.

\bibitem{Hong-MR1050529}
Jooyoo Hong, Yoonbai Kim, and Pong~Youl Pac.
\newblock Multivortex solutions of the abelian {C}hern-{S}imons-{H}iggs theory.
\newblock {\em Phys. Rev. Lett.}, 64(19):2230--2233, 1990.

\bibitem{Jackiw-MR1050530}
R.~Jackiw and Erick~J. Weinberg.
\newblock Self-dual {C}hern-{S}imons vortices.
\newblock {\em Phys. Rev. Lett.}, 64(19):2234--2237, 1990.

\bibitem{K-MR1193342}
Michael K.-H. Kiessling.
\newblock Statistical mechanics of classical particles with logarithmic
  interactions.
\newblock {\em Comm. Pure Appl. Math.}, 46(1):27--56, 1993.

\bibitem{YLi-MR1673972}
Yan~Yan Li.
\newblock Harnack type inequality: the method of moving planes.
\newblock {\em Comm. Math. Phys.}, 200(2):421--444, 1999.

\bibitem{LS-MR1322618}
Yan~Yan Li and Itai Shafrir.
\newblock Blow-up analysis for solutions of {$-\Delta u=Ve^u$} in dimension
  two.
\newblock {\em Indiana Univ. Math. J.}, 43(4):1255--1270, 1994.

\bibitem{LL1-MR2265623}
Chang-Shou Lin and Marcello Lucia.
\newblock Uniqueness of solutions for a mean field equation on torus.
\newblock {\em J. Differential Equations}, 229(1):172--185, 2006.

\bibitem{LL2-MR2352519}
Chang-Shou Lin and Marcello Lucia.
\newblock One-dimensional symmetry of periodic minimizers for a mean field
  equation.
\newblock {\em Ann. Sc. Norm. Super. Pisa Cl. Sci. (5)}, 6(2):269--290, 2007.

\bibitem{Moser-MR0301504}
J.~Moser.
\newblock A sharp form of an inequality by {N}. {T}rudinger.
\newblock {\em Indiana Univ. Math. J.}, 20:1077--1092, 1970/71.

\bibitem{N-Tarantello-MR1664542}
Margherita Nolasco and Gabriella Tarantello.
\newblock On a sharp {S}obolev-type inequality on two-dimensional compact
  manifolds.
\newblock {\em Arch. Ration. Mech. Anal.}, 145(2):161--195, 1998.

\bibitem{Ricciardi-Tarantello-MR1664758}
Tonia Ricciardi and Gabriella Tarantello.
\newblock On a periodic boundary value problem with exponential nonlinearities.
\newblock {\em Differential Integral Equations}, 11(5):745--753, 1998.

\bibitem{Struwe-Tarantello-MR1619043}
Michael Struwe and Gabriella Tarantello.
\newblock On multivortex solutions in {C}hern-{S}imons gauge theory.
\newblock {\em Boll. Unione Mat. Ital. Sez. B Artic. Ric. Mat. (8)},
  1(1):109--121, 1998.

\bibitem{Tarantello-MR1400816}
Gabriella Tarantello.
\newblock Multiple condensate solutions for the {C}hern-{S}imons-{H}iggs
  theory.
\newblock {\em J. Math. Phys.}, 37(8):3769--3796, 1996.

\bibitem{T-MR0216286}
Neil~S. Trudinger.
\newblock On imbeddings into {O}rlicz spaces and some applications.
\newblock {\em J. Math. Mech.}, 17:473--483, 1967.

\end{thebibliography}

\end{document}